\documentclass[12pt]{article}

\usepackage{amsmath,amssymb,amsthm,amsfonts,latexsym,tikz,hyperref}

\newcommand{\wmaj}{\stackrel{\maj}{\equiv}}
\newcommand{\winv}{\stackrel{\inv}{\equiv}}

\DeclareMathOperator{\exc}{exc}
\newcommand{\Ct}{\tilde{C}}
\newcommand{\iotahat}{\hat{\iota}}
\newcommand{\iotacheck}{\check{\iota}}

\newcounter{x}
\newcommand{\permGrid}[1] 
{
\setcounter{x}{0}
\foreach \a in {#1}
{
\addtocounter{x}{1}
\filldraw(\value{x},\a) circle (2pt);
}
\draw (1,1)--(1,\value{x})--(\value{x},\value{x})--(\value{x},1)--(1,1);
}

\newcommand{\permGridInfl}[2] 
{
\newcounter{x}
\setcounter{x}{0}
\newcounter{y}
\setcounter{y}{0}
\foreach \b in {#2}
{
\setcounter{y}{\value{y}+\b}
}
\draw (1,1)--(1,\value{y})--(\value{y},\value{y})--(\value{y},1)--(1,1);
}

\newtheorem{thm}{Theorem}[section]
\newtheorem{prop}[thm]{Proposition}
\newtheorem{cor}[thm]{Corollary}
\newtheorem{lem}[thm]{Lemma}
\newtheorem{conj}[thm]{Conjecture}
\newtheorem{exa}[thm]{Example}
\newtheorem{question}[thm]{Question}

\newcommand{\ben}{\begin{enumerate}}
\newcommand{\een}{\end{enumerate}}
\newcommand{\ble}{\begin{lem}}
\newcommand{\ele}{\end{lem}}
\newcommand{\bth}{\begin{thm}}
\renewcommand{\eth}{\end{thm}}
\newcommand{\bpr}{\begin{prop}}
\newcommand{\epr}{\end{prop}}
\newcommand{\bco}{\begin{cor}}
\newcommand{\eco}{\end{cor}}
\newcommand{\bcon}{\begin{conj}}
\newcommand{\econ}{\end{conj}}
\newcommand{\bde}{\begin{defn}}
\newcommand{\ede}{\end{defn}}
\newcommand{\bex}{\begin{exa}}
\newcommand{\eex}{\end{exa}}
\newcommand{\barr}{\begin{array}}
\newcommand{\earr}{\end{array}}
\newcommand{\btab}{\begin{tabular}}
\newcommand{\etab}{\end{tabular}}
\newcommand{\beq}{\begin{equation}}
\newcommand{\eeq}{\end{equation}}
\newcommand{\bea}{\begin{eqnarray*}}
\newcommand{\eea}{\end{eqnarray*}}
\newcommand{\bce}{\begin{center}}
\newcommand{\ece}{\end{center}}
\newcommand{\bpi}{\begin{picture}}
\newcommand{\epi}{\end{picture}}
\newcommand{\bpp}{\begin{picture}}
\newcommand{\epp}{\end{picture}}
\newcommand{\bfi}{\begin{figure} \begin{center}}
\newcommand{\efi}{\end{center} \end{figure}}
\newcommand{\bprf}{\begin{proof}}
\newcommand{\eprf}{\end{proof}\medskip}

\newcommand{\bsl}{\begin{slide}{}}
\newcommand{\esl}{\end{slide}}
\newcommand{\bfr}{\begin{frame}}
\newcommand{\efr}{\end{frame}}

\newcommand{\pf}{{\bf Proof.\hspace{5pt}}}
\newcommand{\hqed}{\hfill \qed}
\newcommand{\hqedm}{\hfill \qed \medskip}

\newcommand{\eqqed}[1]{$\rule{1ex}{0ex}\hfill{\dil#1}\hfill\qed$}

\newcommand{\hso}[1]{\hspace{-1pt}}
\newcommand{\vs}[1]{\vspace{#1}}
\newcommand{\qmq}[1]{\quad\mbox{#1}\quad}

\newcommand{\emp}{\emptyset}

\newcommand{\sbs}{\subset}
\newcommand{\sbe}{\subseteq}

\newcommand{\spe}{\supseteq}

\newcommand{\case}[4]{\left\{\barr{ll}#1&\mbox{#2}\\#3&\mbox{#4}\earr\right.}

\def\<{\langle}
\def\>{\rangle}

\newcommand{\spn}[1]{\langle{#1}\rangle}
\newcommand{\ree}[1]{(\ref{#1})}

\newcommand{\ra}{\rightarrow}

\newcommand{\be}{\beta}

\newcommand{\de}{\delta}
\newcommand{\ep}{\epsilon}

\newcommand{\io}{\iota}

\newcommand{\la}{\lambda}

\newcommand{\si}{\sigma}
\renewcommand{\th}{\theta}

\newcommand{\bbN}{{\mathbb N}}

\newcommand{\fS}{{\mathfrak S}}

\DeclareMathOperator{\Av}{Av}

\DeclareMathOperator{\des}{des}
\DeclareMathOperator{\Des}{Des}

\DeclareMathOperator{\inv}{inv}
\DeclareMathOperator{\Inv}{Inv}

\DeclareMathOperator{\maj}{maj}

\DeclareMathOperator{\st}{st}

\newcommand{\dil}{\displaystyle}

\begin{document}
\pagestyle{plain}

\title{Permutation patterns and statistics
}
\author{
Theodore Dokos\thanks{Research partially supported by NSA grant H98230-11-1-0222 and by Michigan State University}\\[-5pt]
\small Department of Mathematics, The Ohio State University\\[-5pt]
\small 100 Math Tower, 231 West 18th Avenue\\[-5pt]
\small Columbus, OH 43210-1174, USA, {\tt t.dokos@gmail.com}
\\[5pt]
Tim Dwyer\thanks{Research partially supported by NSA grant H98230-11-1-0222 and by Michigan State University}\\[-5pt]
\small Department of Mathematics,  University of Florida,\\[-5pt]
\small 358 Little Hall, PO Box 118105\\[-5pt]
\small Gainesville, FL 32611-8105, USA, {\tt jtimdwyer@gmail.com}
\\[5pt]
Bryan P. Johnson\thanks{Research partially supported by NSA grant H98230-11-1-0222 and by Michigan State University}\\[-5pt]
\small Department of Mathematics, Michigan State University,\\[-5pt]
\small East Lansing, MI 48824-1027, USA, {\tt john2954@msu.edu}
\\[5pt]
Bruce E. Sagan\thanks{Research partially supported by NSA grant H98230-11-1-0222 and by Michigan State University}\\[-5pt]
\small Department of Mathematics, Michigan State University,\\[-5pt]
\small East Lansing, MI 48824-1027, USA, {\tt sagan@math.msu.edu}
\\[5pt]
Kimberly Selsor\thanks{Research partially supported by NSA grant H98230-11-1-0222 and by Michigan State University}\\[-5pt]
\small Department of Mathematics, University of South Carolina\\[-5pt]
\small LeConte College, 1523 Greene Street\\[-5pt]
\small Columbia, SC 29208, USA, {\tt selsork@email.sc.edu}
}

\date{\today\\[10pt]
	\begin{flushleft}
	                                       \\[5pt]
	\end{flushleft}}

\maketitle

\begin{abstract}
Let $\fS_n$ denote the symmetric group of all permutations of $\{1,2,\ldots,n\}$ and let $\fS=\cup_{n\ge0} \fS_n$.  If $\Pi\sbe\fS$  is a set of permutations, then we let $\Av_n(\Pi)$ be the set of permutations in $\fS_n$ which avoid every permutation of $\Pi$ in the sense of pattern avoidance.  One of the celebrated notions in pattern theory is that of Wilf-equivalence, where $\Pi$ and $\Pi'$ are \emph{Wilf equivalent} if $\#\Av_n(\Pi)=\#\Av_n(\Pi')$ for all $n\ge0$.  In a recent paper, Sagan and Savage proposed studying a $q$-analogue of this concept defined as follows.  Suppose $\st:\fS\ra\{0,1,2,\ldots\}$ is a permutation statistic and consider the corresponding generating function $F_n^{\st}(\Pi;q)=\sum_{\si\in\Av_n(\Pi)} q^{\si}$. 
Call $\Pi,\Pi'$ \emph{st-Wilf equivalent\/} if $F_n^{\st}(\Pi;q)=F_n^{\st}(\Pi';q)$ for all $n\ge0$.  We present the first in-depth study of this concept for the $\inv$ and $\maj$ statistics.  In particular, we determine all $\inv$- and $\maj$-Wilf equivalences for any $\Pi\sbe\fS_3$.  This leads us to consider various $q$-analogues of the Catalan numbers, Fibonacci numbers, triangular numbers, and powers of two.  Our proof techniques use lattice paths, integer partitions, and Foata's second fundamental bijection.  We also answer a question about Mahonian pairs raised in the Sagan-Savage article.
\end{abstract}

\section{Introduction}
\label{i}

We begin with some well-known definitions.
Let $\fS_n$ denote the set of permutations $\pi=a_1 a_2\ldots a_n$ of $[n]=\{1,2,\ldots,n\}$, and let $\fS=\cup_{n\ge0}\fS_n$.  
We denote by $\pi(k)$ the entry  $a_k$.  And if $\pi\in\fS_n$ then we say that $\pi$ has \emph{length $n$}.
 Two sequences of distinct integers, $a_1a_2\ldots a_n$ and $b_1 b_2\ldots b_n$, are said to be \emph{order isomorphic} whenever they satisfy
$$
\mbox{$a_i<a_j$ if and only if $b_i<b_j$}
$$
for all $1\le i<j\le n$.  We say that \emph{$\si\in\fS_n$ contains a copy of $\pi\in\fS_k$ as a pattern} if there is a subsequence of $\si$ order isomporphic to $\pi$.  For example, $\si=436152$ contains the pattern $\pi=132$ because of the copy $365$.  On the other hand, $\si$ \emph{avoids} $\pi$ if it does not contain $\pi$ and we let
$$
\Av_n(\pi)=\{\si\in\fS_n\ |\ \mbox{$\si$ avoids $\pi$}\}.
$$
More generally, if $\Pi\sbe\fS$ then we define $\Av_n(\Pi)=\cap_{\pi\in\Pi} \Av_n(\pi)$ so that these permutations avoid every element of $\Pi$.  Just as with symmetric groups, it will be useful to let
$$
\Av(\Pi)=\bigcup_{n\ge0} \Av_n(\Pi).
$$

We say that $\Pi, \Pi'$ are \emph{Wilf equivalent} and write $\Pi\equiv\Pi'$  if, for all $n\ge0$, 
$$
\#\Av_n(\Pi)=\#\Av_n(\Pi')
$$
where the number sign denotes cardinality.   
Given a specific set $\Pi$, we will  drop the curly brackets enclosing the set when writing $\Av_n(\Pi)$ and in similar notations below.
A famous result in this area states that $\pi\equiv\pi'$ for any $\pi,\pi'\in\fS_3$.

In a recent paper~\cite{ss:mp}, Sagan and Savage proposed studying a $q$-analogue of Wilf equivalence defined as follows.
A \emph{permutation statistic}  is a function
$$
\st:\fS\ra\bbN
$$
where $\bbN$ is the nonnegative integers.   Associated with any statistic $\st$ and any $\Pi\sbe\fS$ we have the generating function
$$
F_n^{\st}(\Pi;q)=\sum_{\si\in\Av_n(\Pi)} q^{\st\si}.
$$
Call $\Pi,\Pi'$ \emph{st-Wilf equivalent}, written $\Pi\stackrel{\st}{\equiv}\Pi'$, if
$$
F_n^{\st}(\Pi;q)=F_n^{\st}(\Pi';q)
$$
for all $n\ge0$.  We denote the st-Wilf equivalence class of $\Pi$ by $[\Pi]_{\st}$. 
Note that st-Wilf equivalence implies Wilf equivalence since setting $q=1$ gives 
$F_n^{\st}(\Pi;1)=\#\Av_n(\Pi)$.  In particular, if $\pi\stackrel{\st}{\equiv}\pi'$ then $\pi$ and $\pi'$ must be in 
the same symmetric group since this is true of Wilf equivalence.

Our goal is to study this concept and related ideas for two famous permutation statistics: the inversion number and the major index.   The set of \emph{inversions} of $\si=a_1a_2\ldots a_n$ is
$$
\Inv \si = \{(i,j)\ |\ \mbox{$i<j$ and $a_i>a_j$}\}.
$$
So $\Inv\si$ records pairs of indices where the corresponding elements in $\si$ are out of order.  The \emph{inversion number} of $\si$ is just
$$
\inv\si=\#\Inv\si.
$$
By way of illustration, if $\si=41523$ 
$$
\Inv\si=\{(1,2),\ (1,4),\ (1,5),\ (3,4),\ (3,5)\}
$$
since $a_1=4>a_2=1$ and so forth, giving $\inv\si = 5$.

To define the major index, we must first define descents.  Permutation $\si=a_1a_2\ldots a_n$  has \emph{descent set} 
$$
\Des\si=\{i\ |\ a_i>a_{i+1}\}
$$ 
and \emph{descent number}
$$
\des\si =\#\Des\si.
$$
Descents keep track of the first index in inversion pairs consisting of adjacent elements.   The \emph{major index} of $\si$ is given by
$$
\maj\si=\sum_{i\in\Des\si} i.
$$ 
Continuing our example from the previous paragraph, we see that 
$$
\Des\si=\{1,3\}
$$
and so $\maj\si=1+3=4$.  

The statistics $\inv$ and $\maj$ are intimately connected and, in fact, are equidistributed over $\fS_n$ in the sense that
$$
\sum_{\si\in\fS_n} q^{\inv\si}=\sum_{\si\in\fS_n} q^{\maj\si}.
$$
A statistic whose distribution over $\fS_n$ equals this one is called \emph{Mahonian}.
Since we will use them so often, we will adopt the following abbreviations for the generating functions for these statistics over avoidance sets
$$
I_n(\Pi;q)=F_n^{\inv}(\Pi;q) \qmq{and} M_n(\Pi;q)=F_n^{\maj}(\Pi;q).
$$

This  article is devoted to a study of $I_n(\Pi)$ and $M_n(\Pi)$ and, by the end, we will have completely classified all these polynomials for all $\Pi\sbe\fS_3$ just as Simion and Schmidt~\cite{ss:rp} did with the corresponding cardinalities when $q=1$.  Along the way we will meet a number of interesting $q$-analogues for known combinatorial quantities such as Catalan numbers, Fibonacci numbers, triangular numbers, and powers of two.  We will establish connections with lattice paths, integer partitions, and Foata's  second fundamental bijection.  We will also answer a question of Sagan and Savage~\cite{ss:mp} concerning certain objects called Mahonian pairs.  The rest of this paper is structured as follows.  In the following section we will consider inv- and maj-Wilf equivalence in the case $\#\Pi=1$.
Section~\ref{1q} will talk about $q$-analogues of the Catalan numbers arising in this context.  In Sections~\ref{2I} and~\ref{2M} we will look at $I_n(\Pi)$ and $M_n(\Pi)$, respectively, when $\Pi\sbs\fS_3$ has cardinality 2.  The next two sections investigate the analogous problem for $\#\Pi=3$.  The final section contains concluding remarks.  Various conjectures are scattered throughout the paper.

\section{Equivalence for single permutations}
\label{1e}

\bfi
\begin{tikzpicture}
\permGrid{1,3,2}
\end{tikzpicture}
\hspace{50pt}
\begin{tikzpicture}
\draw (0,0) -- (2,0) -- (2,2) -- (0,2) -- (0,0);
\draw (2/3,0) -- (2/3,2/3) -- (0,2/3);
\draw (2,2/3) -- (4/3, 2/3) -- (4/3,2);;
\draw (2,4/3) -- (2/3,4/3) -- (2/3,2);
\node at (1/3,1/3){$\si_1$};
\node at (1,5/3){$\si_2$};
\node at (5/3,1){$\si_3$};
\end{tikzpicture}
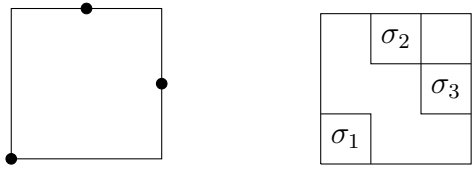
\caption{The diagram of $132$ (left) and $132[\si_1,\si_2,\si_3]$ (right)}
\label{132}
\efi

It will be useful to have a geometric way to describe permutations.  We can diagram $\si=a_1 a_2\ldots a_n \in\fS_n$ using the points $(1,a_1),\ (2,a_2),\ \ldots,\ (n,a_n)$ in the first quadrant of the Cartesian plane.  For example, the diagram of $\si=132$ is in Figure~\ref{132} where we have enclosed the region $1\le x,y\le 3$ in a square.

Inflations of permutations will play an  important role in our work.  Given  permutations $\pi=a_1 a_2\ldots a_k\in\fS_k$ and 
$\si_1,\si_2,\ldots,\si_k\in\fS$, the \emph{inflation} of $\pi$ by the $\si_i$ is the permutation $\pi[\si_1,\si_2,\ldots,\si_k]$ whose diagram is obtained from the diagram of $\pi$ by replacing the point $(i,a_i)$ by the diagram of a permutation order isomorphic to $\si_i$ for $1\le i\le k$. 
The $\si_i$ are called the \emph{components} of the inflation.
 Figure 1 displays the a schematic diagram of $132[\si_1,\si_2,\si_3]$.  As a specific example, we have
$132[21,1,213]=216435$.  Note that we permit $\si_i=\ep$, the empty permutation which has the effect of removing the corresponding point of $\pi$.  By way of illustration $132[\ep,1,213]=4213$.

Since the diagram of $\si$ lies in a square, we can let the dihedral group of the square, $D_4$,  act on permutations.  To set notation, let
$$
D_4=\{R_0,R_{90},R_{180},R_{270},r_{-1},r_0,r_1,r_\infty\}
$$
where $R_\th$ is rotation counterclockwise through an angle of $\th$ degrees and $r_m$ is reflection in a line of slope $m$.  A couple of these rigid motions have easy descriptions in terms of the one line notation for permutations.  If $\si=a_1 a_2\ldots a_n$ then its 
\emph{reversal} is
$$
\si^r = a_n\ldots a_2 a_1 = r_\infty(\si),
$$
and its \emph{complement} is
$$
\si^c = (n+1-a_1), (n+1-a_2),\ldots, (n+1-a_n)=r_0(\si).
$$
We will apply these operations to sets of permutations by doing so element-wise.

We will now investigate what these symmetries do to the inversion polynomial.
Clearly each $f\in D_4$ either preserves inversions and non-inversions, or interchanges the two.  Since the total number of pairs $(i,j)$ with $1\le i<j\le n$ is ${n\choose 2}$, we immediately have the first part of following lemma.  The second part  is obtained from the first by standard generating function manipulations and the fact that $\si$ avoids $\pi$ if and only if $f(\si)$ avoids $f(\pi)$.
\ble
\label{inv:lem}
For any $\si\in\fS$ we have
$$
\inv f(\si)=
\case{\inv\si}{if $f\in\{R_0, R_{180}, r_{-1}, r_1\}$,}
{\dil{n\choose 2}-\inv\si}{if $f\in\{R_{90},R_{270},r_0,r_\infty\}$.\rule{0pt}{20pt}}
$$
It follows that
$$
I_n(f(\pi);q)=
\case{I_n(\pi;q)}{if $f\in\{R_0, R_{180}, r_{-1}, r_1\}$,}
{q^{{n\choose2}}I_n(\pi;q^{-1})}{if $f\in\{R_{90},R_{270},r_0,r_\infty\}$.\rule{0pt}{15pt}}
$$
for all $n\ge0$\hqedm
\ele

We call the $f\in\{R_0, R_{180}, r_{-1}, r_1\}$ \emph{inv preserving} and the $f\in\{R_{90},R_{270},r_0,r_\infty\}$ \emph{inv reversing}.  Translating the previous result into the language of inv-Wilf equivalence gives the following.
\bco
For every permutation $\pi$ and $f\in\{R_0, R_{180}, r_{-1}, r_1\}$ we have
\beq
\label{winv:eq}
f(\pi)\winv \pi.
\eeq
Furthermore, for any $f\in D_4$, we have that $\pi\winv\pi'$ implies $f(\pi)\winv f(\pi')$.\hqedm
\eco

Because of equation~\ree{winv:eq}, we call equivalences arising from these symmetries \emph{trivial}.  Restricting our attention to $\fS_3$, the above corollary gives
$$
132\winv 213 \qmq{and} 231\winv 312.
$$
To see that these are the only inv-Wilf equivalences, note that if $\pi\in\fS_k$ then $\Av_k(\pi)=\fS_k-\{\pi\}$.  It follows that
$$
I_k(\pi;q)=I_k-q^{\inv\pi}
$$
where $I_k=\sum_{\si\in\fS_k} q^{\inv\si}$.  So $\pi\winv\pi'$ forces $\inv\pi=\inv\pi'$.  Since in $\fS_3$ the latter equation only holds for the pairs $132, 213$ and $231, 312$, we have proved the following result.
\bth
\label{inv3:th}
We have
\bea
\left[123\right]_{\inv} &=& \{123\},\\
\left[321\right]_{\inv} &=& \{321\},\\
\left[132\right]_{\inv}  &=& \{132,213\},\\
\left[231\right]_{\inv} &=& \{231,312\},
\eea
and these are all the inv-Wilf equivalence classes involving single elements of $\fS_3$.\hqedm
\eth

In fact, we believe that all inv-Wilf equivalences are trivial.  The following conjecture has been verified by computer for single permutations up through $\fS_5$  and for all  subsets of permutations of $\fS_3$.  Note that if $\Pi=\{\pi_1,\ldots,\pi_l\}$ then we let $f(\Pi)=\{f(\pi_1),\ldots,f(\pi_l)\}$.  In a similar vein, we will let other operators distribute over sets, e.g., if $S=\{a_1,\ldots,a_l\}$ is a set of integers then $n+S=\{n+a_1,\ldots,n+a_l\}$.
\bcon
For $\Pi,\Pi'\sbe\fS_n$ we have $\Pi\winv\Pi'$ if and only if $f(\Pi)=\Pi'$ for some $f\in\{R_0, R_{180}, r_{-1}, r_1\}$.
\econ

We now turn our attention to the maj polynomial.  Here the situation is more complicated because there are no maj-preserving symmetries in $D_4$.   In fact, the only symmetry which permits one to deduce something about maj-Wilf equivalence is complementation.  
The next result is proved in much the same way as Lemma~\ref{inv:lem} and so its demonstration is omitted.
\ble
\label{maj:lem}
For any $\si\in\fS$ we have
$$
\maj\si^c={n\choose2} - \maj\si.
$$
It follows that
$$
M_n(\pi^c;q)=q^{{n\choose2}} M_n(\pi,q^{-1})
$$
for all $n\ge0$.  So $\pi\wmaj\pi'$ if and only if $\pi^c\wmaj(\pi')^c$.  \hqedm
\ele

We now come to our first non-trivial equivalences.  In the proof as well as later we will say that a map $f$ on permutations is \emph{Des preserving} if $\Des\si=\Des f(\si)$ for all $\si$ in the domain of $f$.  Note that such a map must preserve both the $\des$ and $\maj$ statistics in the same sense.
\bth
\label{maj3:th}
We have
\bea
\left[123\right]_{\maj} &=& \{123\},\\
\left[321\right]_{\maj} &=& \{321\},\\
\left[132\right]_{\maj}  &=& \{132,231\},\\
\left[213\right]_{\maj} &=& \{213,312\},
\eea
and these are all the maj-Wilf equivalence classes involving single elements of $\fS_3$.\hqedm
\eth
\pf  Using arguments like those before Theorem~\ref{inv3:th}, we see that $132\wmaj 231$ and $213\wmaj 312$ are the only possible maj-Wilf equivalences in $\fS_3$.  Also, the previous lemma shows that it suffices to prove one of these two equivalences, and we will choose the former.

\bfi
\begin{tikzpicture}
\draw (0,0) -- (4,0) -- (4,4) -- (0,4) -- (0,0);
\draw (2.4,0) -- (2.4,2) -- (4,2);
\draw (0,2) -- (1.6,2) -- (1.6,4);
\node at (.8,3){$\si_1$};
\node at (3.2,1){$\si_2$};
\filldraw(2,4) circle(3pt);
\node at (2,4.4){$n$};
\end{tikzpicture}
\hspace{50pt}
\begin{tikzpicture}
\draw (0,0) -- (4,0) -- (4,4) -- (0,4) -- (0,0);
\draw (2.4,4) -- (2.4,2) -- (4,2);
\draw (0,2) -- (1.6,2) -- (1.6,0);
\node at (.8,1){$f(\si_1)$};
\node at (3.2,3){$f(\si_2)$};
\filldraw(2,4) circle(3pt);
\node at (2,4.4){$n$};
\end{tikzpicture}
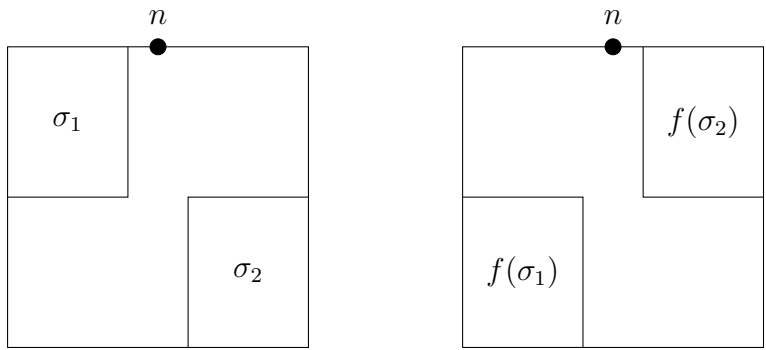
\caption{An element of $Av_n(132)$  (left) and its image under $f$ (right)}
\label{Av(132)}
\efi

We first note that $\si\in\Av_n(132)$ if and only if $\si=231[\si_1,1,\si_2]$ for some $\si_1,\si_2\in\Av(132)$ as shown in the diagram on the left in Figure~\ref{Av(132)}.  This is well known, but since we will be using such characterizations frequently in the sequel, we will present a proof this time and leave future demonstrations to the reader.

Suppose first that $\si\in\Av_n(132)$ and write $\si=\si_L n\si_R$ where $\si_L$ and $\si_R$ are the subsequences of $\si$ to the left and right of $n$, respectively.  To avoid having a copy of $132$ where $n$ acts as the $3$, we must have every element of $\si_L$ larger than every element of $\si_R$.  It follows that $\si=231[\si_1,1,\si_2]$ where $\si_1$ is order isomorphic to $\si_L$ and similarly with $\si_2$ and $\si_R$.  Clearly $\si_1$ and $\si_2$ must avoid $132$ since $\si$ does, and this completes the forward direction of the characterization.

Now suppose $\si$ is an inflation of the given form and that, towards a contradiction, it contains a copy of $132$.  Consider the element $a$ of $\si$ playing the role of the $2$.  If $a$ comes from $\si_1$ then the whole copy is in $\si_1$, contradicting the fact that $\si_1\in\Av(132)$.  It is impossible for $n$ to be the $2$.  So the only other choice is for $a$ to come from $\si_2$.  But since $n$ as well as every element of $\si_1$ is larger then every element of $\si_2$, the element playing the role of the $1$ in the $132$ copy must also come from $\si_2$, placing the whole copy in $\si_2$.  This contradicts $\si_2\in\Av(132)$ and finishes the characterization.  Note that exactly the same line of reasoning shows that $\si\in\Av_n(231)$ if and only if $\si=132[\si_1,1,\si_2]$ where $\si_1,\si_2\in\Av(231)$.

We will now prove $132\wmaj 231$ by inductively constructing a $\Des$-preserving bijection $f:\Av_n(132)\ra\Av_n(231)$ for all $n\ge0$.  For $n=0$ we define $f(\ep)=\ep$.  Suppose that $f$ has been defined for permutations with fewer than $n$ elements and consider $\si\in\Av_n(132)$.  If $\si=231[\si_1,1,\si_2]$ 
then define the image of $\si$ under $f$ to be 
$$
f(\si)=132[f(\si_1),1,f(\si_2)].
$$
This is illustrated by the right-hand diagram  in Figure~\ref{Av(132)}.  By induction and the characterization  above, we have that $f(\si)\in\Av_n(231)$.  So $f$ does indeed map the domain into the desired range.  It is also clear inductively that $f$ preserves the descent set and so preserves maj.
It is an easy matter to construct an inverse for $f$, thus the map is bijective and we are done.
\hqedm

The two maj-Wilf equivalences proved in the previous theorem seem to be part of a larger picture.  To state this result, we use the notation 
$$
\iota_n=12\ldots n
$$ 
and
$$
\de_n=n\ldots21
$$ 
for the monotone increasing and decreasing permutations.  We will suppress the subscript if we do not wish to specify the number of elements in the permutation.
\bcon
For all $m,n\ge0$  we have
$$
132[\iota_m,1,\de_n]\wmaj 231[\iota_m,1,\de_n].
$$
Equivalently
$$
213[\de_m,1,\iota_n]\wmaj 312[\de_m,1,\iota_n].
$$
Also, for $n\ge5$ these are the only maj-Wilf equivalences of the form $\pi\wmaj \pi'$ with $\pi\neq\pi'$.
\econ

The maj-Wilf equivalences in this conjecture have been checked by computer  up to $\fS_7$.  We have also checked that there are no other such equivalences in $\fS_5$.  In $\fS_4$, however, there appear to be some other sporadic ones.
Aside from the two sets of equivalences below, all other maj-Wilf equivalence classes in $\fS_4$ are singletons.
\bcon
We have
$$
[1423]_{\maj}=\{1423, 2314, 2413\}.
$$
Equivalently
$$
[3142]_{\maj}=\{3142, 3241, 4132\}.
$$
\econ

\section{$q$-Catalan numbers}
\label{1q}

The \emph{Catalan numbers} can be defined as
$$
C_n=\frac{1}{n+1}{2n\choose n}.
$$
It is well known that $\#\Av_n(\pi)=C_n$ for any $\pi\in\fS_3$.  So the polynomials $I_n(\pi;q)$ and $M_n(\pi;q)$ for $\pi\in\fS_3$ are $q$-analogues of these numbers in that they evaluate to $C_n$ when setting $q=1$.
In view of Lemmas~\ref{inv:lem} and~\ref{maj:lem} as well as Theorems~\ref{inv3:th} and~\ref{maj3:th}, only four of these polynomials are really different  up to a change of variable.  So we will concentrate on the cases when $\pi=312$ and $321$.

The Catalan numbers can also be defined as the unique solution to the initial condition $C_0=1$ and recursion
$$
C_n=\sum_{k=0}^{n-1} C_k C_{n-k-1}
$$
for  $n\ge1$.  Cartlitz and Riordan~\cite{cr:tel} defined a $q$-analogue $C_n(q)$ of $C_n$.  For our purposes, it will be easier to work with the polynomials 
\beq
\label{CtC}
\Ct_n(q)=q^{{n\choose2}} C_n(q^{-1})
\eeq 
gotten by reversing the order of the coefficients of $C_n(q)$.  These can be defined by $\Ct_0(q)=1$ and, for $n\ge1$,
\beq
\label{Ct}
\Ct_n(q)=\sum_{k=0}^{n-1} q^k \Ct_k(q) \Ct_{n-k-1}(q).
\eeq
The Carlitz-Riordan polynomials have been studied by numerous authors.  They have a close connection with the space of diagonal harmonics and Haglund's text~\cite{hag:qtc} gives more information and references. 

It turns out that some of our polynomials are related to those of  Carlitz and Riordan.
\bth
We have $I_0(312;q)=1$  and
\beq
\label{I(312)}
I_n(312;q)=\sum_{k=0}^{n-1} q^k I_k(312;q) I_{n-k-1}(312;q).
\eeq
for $n\ge1$.  It follows that, for $n\ge0$,
$$
\barr{l}
C_n(q) = I_n(132;q)=I_n(213;q),\\
\Ct_n(q)=I_n(231;q)=I_n(312;q).
\earr
$$
\eth
\pf
To prove the recursion, suppose that $\si\in\Av_n(312)$ so that $\si=213[\si_1,1,\si_2]$.  If $\si(k+1)=1$ for some $0\le k<n$
then $\si_1\in\Av_k(312)$, $\si_2\in\Av_{n-k-1}(312)$, and
$$
\inv\si=k+\inv\si_1+\inv\si_2.
$$
Rewriting the previous equation in terms of generating functions gives the desired  result. 

The statements about $\Ct_n(q)$   now follow by comparing~\ree{Ct} and~\ree{I(312)} and applying Theorem~\ref{inv3:th}.
Those about $C_n(q)$ are then obtained from~\ree{CtC} and Lemma~\ref{inv:lem}.
\hqedm

Although a similar recursion seems to hold for $I_n(321)$, we have been unable to prove it and so state it as a conjecture.
\bcon
\label{I321:con}
We have, for all $n\ge1$,
$$
I_n(321;q)=I_{n-1}(321;q)+\sum_{k=0}^{n-2} q^{k+1} I_k(321;q) I_{n-k-1}(321;q).
$$
\econ

We have not been able to  discover a recursion for the polynomials $M_n(312;q)$ and it would be interesting to do so.  However, one can get a result by using a bivariate generating function.  For $\Pi\sbe\fS$ we let
$$
M(\Pi;q,t)=\sum_{\si\in\Av_n(\Pi)}q^{\maj\si} t^{\des\si}.
$$
\bth
We have, for $n\ge1$,
$$
M_n(312;q,t)=M_{n-1}(312;q,qt)+\sum_{k=1}^{n-1} q^k t M_k(312;q,t) M_{n-k-1}(312;q,q^{k+1}t).
$$
\eth
\pf
Keeping the notation of the previous proof, we have
$$
\des\si=
\case{\des\si_2}{if $k=0$,}{1+\des\si_1+\des\si_2}{if $k\ge1$,}
$$
and, for all $k\ge0$,
$$
\maj\si = k+\maj\si_1+\maj\si_2+(k+1)\des\si_2.
$$
Turning these statements into an equation for the generating functions completes the proof.
\hqedm

Computer experimentation has yielded no nice recursion for $M_n(321;q,t)$ so we leave this as an open problem
\begin{question}
\label{M321:ques}
Do the $M_n(321;q,t)$ satisfy a $(q,t)$-analogue of the Catalan recursion?
\end{question}

We will now show how one can use one of our $q$-analogues to give a new proof of an interesting result about Catalan numbers.  It is well known that
\beq
\label{parity}
\mbox{$C_n$ is odd if and only if $n=2^k-1$ for some $k\ge0$.}
\eeq
In fact, one can characterize the $2$-adic valuation (highest power of 2) dividing $C_n$.  For a statement and combinatorial proof of this stronger result, see the paper of Deutsch and Sagan~\cite{ds:ccm}.  We will prove a refinement of the ``if" direction of~\ree{parity}.  To do so  it will be convenient to have the notation that, for any polynomial $f(q)$,
$$
\spn{q^i}f(q)=\mbox{the coefficient of $q^i$ in $f(q)$}.
$$
\bth
For all $k\ge0$ we have
$$
\spn{q^i} I_{2^k-1}(321;q)=
\case{1}{if $i=0$,}{\mbox{an even number}}{if $i\ge1$.}
$$
\eth
\pf
The constant coefficient of $I_n(321;q)$ is 1 because $\iota_n$ is the only permutation in $\fS_n$ with no inversions.
And if $n=2^0-1=0$ then all other coefficients are zero, which is even.  So we will proceed by induction on $k$.

Assume the result is true for indices less than $k$ and consider
$$
S_i=\{\si\in\Av_{2^k-1}(321)\ |\ \inv\si=i\}.
$$
For $i\ge1$, we will partition this set into blocks (subsets) each of even size and that will prove the theorem.  Since 
$R_{180}(321)=321$ and $R_{180}$ preserves inv by Lemma~\ref{inv:lem}, the involution $R_{180}$ acts on $S_i$.  

\bfi
\begin{tikzpicture}
\draw (0,0) -- (4,0) -- (4,4) -- (0,4) -- (0,0);
\draw (1.6,0) -- (1.6,1.6) -- (0, 1.6);
\draw (2.4,0) -- (2.4,1.6) -- (4,1.6);
\draw (2.4,4) -- (2.4,2.4) -- (4,2.4);
\draw (0,2.4) -- (1.6,2.4) -- (1.6,4);
\node at (3.2,3.2){$\si_{I}$};
\node at (.8,3.2){$\si_{II}=\emp$};
\node at (.8,.8){$\si_{III}$};
\node at (3.2,.8){$\si_{IV}=\emp$};
\filldraw(2,2) circle(3pt);
\node at (1.5,2){$P$};
\end{tikzpicture}
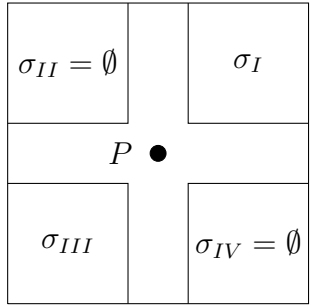
\caption{An element of $\Av_{2^k-1}(321)$ fixed by $R_{180}$}
\label{R_180}
\efi

Suppose $\si\in S_i$ is such that $R_{180}(\si)\neq\si$.  Then let the block containing $\si$ be $B=\{\si,R_{180}(\si)\}$.  By construction, all these blocks have size two which is even.  Now put all of the remaining $\si$ in one block
$$
B'=\{\si\in S_i\ |\ R_{180}(\si)=\si\}.
$$

For $R_{180}(\si)=\si$ to hold, one of the points of the diagram of $\si$ must be the center of rotation $P=(2^{k-1},2^{k-1})$ which splits the diagram up into four parts as in Firgure~\ref{R_180}.  Note that since $\si_{IV}=R_{180}(\si_{II})$ both parts must be empty or nonempty.  If they are both nonempty then this will result in a copy of $321$, so they are both empty.  Also  $\si_{III}=R_{180}(\si_I)$ and so once $\si_I$ is chosen, $\si_{III}$ is determined.  Thus we can characterize the $\si\in B'$ as those of the form $\si=123[R_{180}(\tau),1,\tau]$ where $\tau\in\Av_{2^{k-1}-1}(321)$.  Furthermore $i=\inv\si=2\inv\tau$.  It follows that
$$
\# B' = \spn{q^{i/2}}I_{2^{k-1}-1}(321;q)
$$
and, by induction, this is an even number.  This finishes the proof.
\hqedm

We conjecture that a similar statment can be made for one of the maj polynomial sequences.
\bcon
\label{M321parity}
For all $k\ge0$ we have
$$
\spn{q^i} M_{2^k-1}(321;q)=
\case{1}{if $i=0$,}{\mbox{an even number}}{if $i\ge1$.}
$$
\econ

\
\section{The inversion polynomial for doubletons}
\label{2I}

By the Erd\H{o}s-Szekeres theorem on monotonic subsequences \cite{es:cpg}, any permutation of length 5 or greater will contain an increasing or decreasing subsequence of length 3.  Due to this fact, we have that $I_n(\Pi;q)=0$ (similarly $M_n(\Pi;q)=0$) for $n\geq5$ and $\Pi\spe\{123, 321\}$.  For such $\Pi$ different techniques will be used and so we will not consider the set $\{123,321\}$ here.  Such $\Pi$ will be discussed  in Section \ref{op}.

For $\Pi \sbe \fS$, we define 
\begin{equation*}
\inv \Pi = \{\{\inv \pi | \pi \in \Pi\}\}\\
\end{equation*}
to be the multiset of $\inv \pi$ for each $\Pi$. 
We define $\maj \Pi$ similarly.  For example, the set $\Pi=\{123, 132, 213, 231, 321\}$ has $\inv\Pi = \{\{0, 1, 1, 2, 3\}\}$ and $\maj\Pi = \{\{0, 1, 2, 2, 3\}\}$.
Note that if $\Pi, \Pi' \sbe\fS_n$ for a fixed $n$ and $\Pi\stackrel{\inv}{\equiv}\Pi'$, then $\inv \Pi = \inv \Pi'$ and the same result holds true for $\maj\Pi$.   The argument is similar to the one given before Theorem~\ref{inv3:th} in the case $\#\Pi=\#\Pi'=1$.

In the next result, the inv-Wilf equivalences are obtained from
equation~\ree{winv:eq}.  On the other hand, it is easy to check that $\inv\Pi\neq\inv\Pi'$ for $\Pi$ and $\Pi'$ listed in different classes and so, by the discussion above, there can be no other such equivalences among these doubletons.
\bth
\label{inv2:list}

We have
\bea
\left[123, 132\right]_{\inv} &=& \big\{\{123, 132\}, \{123, 213\}\big\},\\
\left[312, 321\right]_{\inv} &=& \big\{\{231, 321\}, \{312, 321\}\big\},\\
\left[123, 231\right]_{\inv} &=& \big\{\{123, 231\}, \{123, 312\}\big\},\\
\left[213, 321\right]_{\inv} &=& \big\{\{132, 321\}, \{213, 321\}\big\},\\
\left[132, 231\right]_{\inv} &=& \big\{\{132, 231\}, \{132, 312\}, \{213, 231\}, \{213, 312\}\big\}.
\eea
All other inv-Wilf equivalence classes for $\Pi\subset \fS_3$ with $|\Pi|=2$ and $\Pi\neq\{123,321\}$ consist of a single pair.
\hqedm
\eth

The equivalence classes above have been paired, with the class of $\Pi$ appearing in a line adjacent with the class of its complement.  By Lemma~\ref{inv:lem}, $I_n(\Pi;q)$ for one class determines the polynomial for its pair.  So we only compute the inversion polynomial for one doubleton from each pair consisting of a class and it's complement; the reader can easily do the computations for the complementary class.

Before giving the next result, we first must make some remarks about integer partitions. Let $\lambda=(\lambda_1,\lambda_2,\dots,\lambda_k)$ be an \emph{integer partition}, that is, a weakly decreasing sequence of positive integers called \emph{parts}.
We use the notation
$$
|\lambda|=\lambda_1+\lambda_2+\cdots +\lambda_k
$$ 
and 
$$
l(\lambda)= \mbox{the number of parts of $\la$.}
$$ 
If we further suppose that $n > \lambda_1 > \lambda_2 > \cdots > \lambda_k,$ then $\lambda$ is a partition of $|\lambda|$ into \emph{distinct} parts, with no part larger than $n-1$.  We will use the notation $P_{d}(n-1)$ for this set. It is well known that the generating function for these partitions is expressible as a simple product
\beq
\label{part:eq}
\sum_{\lambda\in P_d(n-1)}q^{|\lambda|}t^{l(\la)}=(1+qt)(1+q^2t)\cdots(1+q^{n-1}t).
\eeq
We also define the following permutations which will be useful here  as well as in Section \ref{2M}
\beq
\iotacheck_n =n1\cdots (n-1) = 21[1,\iota_{n-1}], \label{iotacheck}
\eeq
and
\beq
\iotahat_n = 2\cdots n1 = 21[\iota_{n-1},1]. \label{iotahat}
\eeq
Note that the polynomials in the following theorem give two different $q$-analogues of the number $2^{n-1}$.

\bpr
\label{inv2:form}
We have
\begin{align*}
I_n(231, 321;q) &= (1+q)^{n-1}\\
I_n(132, 231;q) &= (1+q)(1+q^{2})\cdots(1+q^{n-1})
\end{align*}
for all $n\geq 1$.
\epr

\pf For the first equation, it is easy to see that $\si$ is determined by its left-right maxima
$m_1<m_2<\cdots<m_j=n$.   In fact, $\sigma=\iota_{j}[\iotacheck_{n_1},\iotacheck_{n_2},\dots,\iotacheck_{n_j}]$ where $n_i=m_i-m_{i-1}$ for $1\le i<k$. (By convention $m_0=0$.)   Note that $\inv\si=n-j$ since each element which is not a left-right maximum participates in exactly one inversion.  And there are ${n-1\choose j-1}$ ways to chose the $m_i$ since $m_j=n$ is always a left-right maximum.  So we have 
$$
I_n(231, 321;q)=\sum_{j=1}^{n}{{n-1}\choose{j-1}}q^{n-j}
=\sum_{j=0}^{n-1}{{n-1}\choose{j}}q^{j}
=(1+q)^{n-1}.
$$

The second equation comes from noticing that $\sigma\in\Av_n(132, 231)$ if and only if $\sigma=\sigma_L1\sigma_R$
where $\sigma_L$ is a decreasing sequence and $\sigma_R$ is an increasing sequence.  We now define a bijection $f:\Av_n(132, 231) \longrightarrow P_{d}(n-1)$.  From the discussion before the theorem, it suffices to show that for $\lambda=f(\sigma)$, we have that $\inv\sigma=|\lambda|$.  Note that  for any entry $i$ in $\sigma_L$, all entries $1, \dots, i-1$ occur after $i$ and that $\sigma_R$ contributes nothing to $\inv\sigma$.  So $\inv\sigma=(\sigma(1)-1)+\dots+(\sigma(k)-1)$.  Defining $f(\sigma)=(\sigma(1)-1, \dots, \sigma(k)-1)$, it is easy to check that $f$ has the desired properties.
\hqedm

For the remaining two equivalence classes, $\left[132, 321\right]_{\inv}$ and $\left[132, 213\right]_{\inv}$, we have not been able to obtain closed forms for their inversion polynomials.  Therefore, only summation formulae will be proved here.

\bpr
\label{inv2:rec}
For all n $\ge$ 1, we have the following recurrences
\begin{align*}
I_n(132, 321;q) &=1 + \sum_{k=1}^{n-1} \frac{q^{k(n-k+1)}-q^k}{q^k-1},\\
I_n(132, 213;q) &= \sum_{k=1}^{n}{q^{k(n-k)}I_{n-k}(132, 213;q)}
\end{align*}
where the latter equation is a recurrence for one of the singletons not listed in $Theorem~\ref{inv2:list}$.
\epr

\pf By considering the placement of $n$, one sees that  $\sigma\in\Av_n(132, 321)$ if and only if $\sigma = 12[\sigma', 1]$ with 
$\sigma'\in\Av_{n-1}(132, 321)$ or $\sigma = 21[\iota_{k}, \iota_{n-k}]$ for some $k$, $0<k<n$.  These two cases result in $\inv\sigma = \inv\sigma'$ or $\inv\sigma = k(n-k)$, respectively.  It follows that, for $n\ge1$,
$I_n(132, 321;q) = I_{n-1}(132, 321;q)+\sum_{k=1}^{n-1}{q^{k(n-k)}}$.  Iterating this recursion gives the desired sum.

For the second equality, note that $\sigma\in\Av_n(132, 213)$ if and only if $\sigma = 21[\iota_k, \sigma']$ where $\sigma'\in\Av_{n-k}(132, 213)$ for some $k\in [n]$.  Therefore, $\inv\sigma = \inv\sigma' +k(n-k)$ and the second recurrence follows.
\hqedm

\section{The major index polynomial for doubletons}
\label{2M}
As discussed in the previous section, we will  consider the set $\{123,321\}$ and its supersets at the end of the paper. For the remaining 14 pairs, there is only one maj-Wilf equivalence class which does not consist of a single pair.
\bth
\label{maj2:list}
We have
$$
\left[ 132,213\right]_{\maj} = \{\{132,213\},\{132,312\},\{213,231\},\{231,312\}\}.
$$
All other maj-Wilf equivalence classes for $\Pi\subset\fS_3$ with $|\Pi|=2$ and $\Pi\neq\{123,321\}$ contain a single pair.
\eth
\pf
The non-equivalences follow from computer-generated example polynomials.  The equivalence of the sets $\{132,213\}, \{132,312\},$ and $\{213,231\}$ will follow from the next proposition, where we will explicitly compute the corresponding polynomials and see that they are the same.  By complementing the equivalence $\{132,213\}\wmaj\{132,312\}$ we obtain  $\{231,312\}\wmaj\{132,312\}$, completing the proof.
\hqedm

We will now explicitly calculate $M_n(\Pi;q,t)$ for the various $\Pi$ under consideration.  For every pair $\Pi$, $\Pi^c$, we will only compute one of the polynomials since the other can be obtained using the fact that $M_n(\Pi^c;q)=q^{\binom{n}{2}}M_n(\Pi;q^{-1})$ by Lemma~\ref{maj:lem}.

\bpr
\label{maj:132,213}
For $n\ge0$ and $\Pi=\{132,213\}$, $\{132,312\}$, or $\{213,231\}$ we have
$$M_n(\Pi;q,t)=(1+qt)(1+q^2t)\cdots(1+q^{n-1}t).$$
\epr
\pf
When $\Pi=\{132,213\}$, in view of~\ree{part:eq}, it suffices to establish a bijection $f:\Av_n(132,213)\ra P_d(n-1)$  such that if $\sigma\in\Av_n(132,212)$ and  $\lambda=f(\sigma)$, then  $\des\sigma=l(\lambda)$ and $\maj \sigma=|\lambda|$.
Recursively applying the characterization of $\Av_n(132,213)$ given in the proof of Proposition~\ref{inv2:rec},  we see that $\sigma$ avoids $132, 213$  if and only if  $\sigma=\delta_k[\iota_{m_1},\iota_{m_2},\dots ,\iota_{m_k}]$.
Thus $\sigma$ is uniquely determined by the elements of its descent set, which we shall denote $\Des\si=\{\lambda_1 > \lambda_2 >\cdots >\lambda_k \}$. But by letting $\lambda=(\lambda_1 , \lambda_2 ,\dots ,\lambda_k )$, we have our bijection which clearly satisfies $l(\lambda)=\des\sigma$ and $|\lambda|=\maj\sigma$.

For $\Pi=\{132,312\}$, it now suffices to show that there is a Des-preserving bijection  
$f:\Av_n(132,312)\ra\Av_n(132,213)$.
Since avoiding $132$ and $312$ forces the permutation to end in $1$ or $n$, we have $\sigma\in\Av_n(132,312)$ if and only if 
$$
\si=21[\si',1]\mbox{ or }12[\si',1]
$$
where $\si'\in\Av_{n-1}(132,312)$.  By induction, we can assume 
$f(\si')=\tau'\in\Av_{n-1}(132,213)$ where the image permutation has the form 
$\tau'=\delta_k[\iota_{m_1},\iota_{m_2},\dots ,\iota_{m_k}]$.  So define
$$
f(\si)=\delta_{k+1}[\iota_{m_1},\iota_{m_2},\dots ,\iota_{m_k},1]\mbox{ or }
\delta_k[\iota_{m_1},\iota_{m_2},\dots ,\iota_{m_k+1}]
$$
where the two cases for $f(\si)$ correspond respectively to the two cases for $\si$ listed above.
It is easy to see that this process can be reversed to obtain an inverse, which proves that $f$ is a bijection. Also $\Des\si'=\Des\tau'$ by induction and  the final ascent or descent  has been preserved by construction.  So $f$ is Des preserving as desired.

For $\Pi=\{213,231\}$, we again wish to construct  a $\Des$-preserving bijection, this time from $\Av_n(213,231)$ to $\Av_n(132,213)$. 
But $\si\in\Av_n(213,231)$ if and only if $\si(1)=1$ or $\si(1)=n$
and $\sigma'=\si(2)\ldots\si(n)\in\Av_{n-1}(213,231)$.  It should now be clear how  the proof of the previous paragraph can be adapted to provide the desired bijection.
\hqedm

The next result is easy to show using the previous characterizations of $\Av_n(132,231)$  and $\Av_n(132,321)$ in the demonstrations of Propositions~\ref{inv2:form} and~\ref{inv2:rec}.  So we omit the proof.  In the result for the second pair, we use the standard
$q$-analogue of a natural number $n$:
$$[n]_q=1+q+\cdots+q^{n-1}.$$

\bpr
\label{maj:132,231}
 For $n\geq 0$, we have
$$
M_n(132,231;q,t)=\sum_{k=0}^{n-1} \binom{n-1}{k}q^{\binom{k+1}{2}}t^{k}
$$
and

\eqqed{
M_n(132,321;q,t)=1+qt\frac{n-[n]_q}{1-q}.
}
\epr

\bpr
\label{maj:213,321}
For $n\geq 0$, we have
$$M_n(213,321;q,t)=1+ t\sum_{k=1}^{n-1} kq^k .$$
\epr
\pf
The permutations in $\Av_n(213,321)$ are those of the form $\si=132[\io_{i_1},\io_{i_2},\io_{i_3}]$.
As usual, let $k=\sigma^{-1}(n)$ and $\si=\si_L n\si_R$.
If $k<n$, then the number of ways of partitioning $\si_L$ is equal to $k$. If $k=n$, then $\sigma=\iota_n$ which contributes the 1 present in our expression. Translating these facts into a generating function completes the proof.\hqedm

It is known that $\#\Av_n(213,321)=1+\binom{n}{2}$, from~\cite{ss:rp}, and thus our polynomial can be interpreted as a $q$-analogue of the expression 
$$\binom{n}{2}=\sum_{k=1}^{n-1}k.$$

The final two sets of permutations do not have known closed forms for their major index polynomials. However, their generating functions 
$$M(\Pi;x)=\sum_{n\geq 0} M_n(\Pi;q,t)x^n$$
are equal to generating functions over binary words. It will be informative to look at these words and their intimate connection with integer partitions.

\begin{figure}
\begin{center}
\begin{tikzpicture}
\filldraw [color=black!10] (0,5) rectangle (2,3) ;
\draw [very thick] (0,0)--(0,2)--(2,2)--(2,3)--(3,3)--(3,5)--(4,5);
\filldraw (0,0) circle (2pt);
\filldraw (0,1) circle (2pt);
\filldraw (0,2) circle (2pt);
\filldraw (1,2) circle (2pt);
\filldraw (2,2) circle (2pt);
\filldraw (2,3) circle (2pt);
\filldraw (3,3) circle (2pt);
\filldraw (3,4) circle (2pt);
\filldraw (3,5) circle (2pt);
\filldraw (4,5) circle (2pt);

\draw (0,2)--(0,5)--(3,5);
\draw (1,2)--(1,5);
\draw (2,3)--(2,5);
\draw (0,4)--(3,4);
\draw (0,3)--(2,3);
\end{tikzpicture}
\end{center}
\caption{A lattice path with Ferrers diagram and Durfee square shown.}
\label{lpath}
\end{figure}
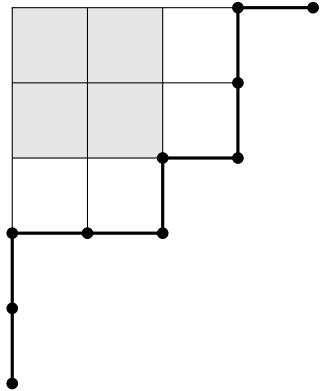

Given a 0-1 word $w$ with $m$ zeros and $n$ ones, we associate with it a \emph{lattice path}, $P(w)$, from $(0,0)$ to $(n,m)$ with each zero representing a unit step to the north, and each one a unit step to the east. In Figure~\ref{lpath} we have the path $001101001$. The unit boxes northwest of the path and contained in the rectangle with diagonal corners $(0,0)$ and $(n,m)$ form the \emph{Ferrers diagram} of the path. This diagram can be interpreted as a partition $\lambda(w)$, with each row of boxes being a part. For example, the partition in Figure~\ref{lpath} is (3,3,2). It is well known~\cite{ss:mp} that 
\beq
\label{invtosumparts}
\inv w=|\lambda(w)|.
\eeq

For a given Ferrers diagram of a path $P(w)$, the \emph{Durfee square} $D(w)$ is the largest square which fits into the partition $\lambda=\lambda(w)$ so that its northwest corner coincides with the northwest corner of the diagram. Its side length is denoted $d(w)$. In Figure~\ref{lpath} the Durfee square is shaded, and $d(w)=2$. 

Our generating function sums will be obtained by decomposing a Ferrers diagram into its Durfee square, the part of the diagram below  the square, and the part of the diagram to the right of the square.  The part below corresponds to  a partition $\be(w)$ where each part has size at most $d(w)$. Additionally, we allow parts of size 0 in order to incorporate any steps to the north that precede all eastern steps. In Figure~\ref{lpath}, we have $\be(w)=(2,0,0)$.   Similarly, if we view the columns to the right of $D(w)$ as parts instead of the rows, we obtain another partition $\rho(w)$.  Referring to Figure~\ref{lpath} again, we see that $\rho(w)=(2,0)$.

Finally, we will use Foata's second fundamental bijection $\phi:\mathbb{N}^*\to\mathbb{N}^*$, where $\mathbb{N}^*$ is the Kleene closure of all finite words formed by elements of $\mathbb{N}$. Although we do not define this map here, the relevant properties will be stated in the next three lemmas. We shall use $v$ to denote a word in the domain of $\phi$, and $w$ a word in the range. Additionally, the length (number of elements) of any word $u\in\mathbb{N}^*$ will be denoted $|u|$. Foata's bijection originally appeared in~\cite{foa:nin} and was constructed to have the following properties.

\ble[\cite{foa:nin}]
\label{thm4.3:foa:nin}
For $v\in\mathbb{N}^*$, we have that 
$$|v|=|\phi(v)|,
$$
and

\vs{5pt}

\eqqed{\inv\phi(v)=\maj v.}
\ele

The next two lemmas are from Sagan and Savage's paper~\cite{ss:mp}. The form for $v$ given in Lemma~\ref{prop2.2:ss:mp} is the general form of a 0-1 word with $\des v=k$. Readers unfamiliar with $\phi$ may take the Lemma as the definition of $\phi$ on $\{0,1\}^*$:

\ble[\cite{ss:mp}]
\label{prop2.2:ss:mp}
Let $v=0^{m_0}1^{n_0}0^{m_1}1^{n_1}\dots 0^{m_k}1^{n_k}\in \{0,1\}^*$, where $m_0,n_k\geq 0$ and $m_i,n_j>0$ for all other $i$ and $j$. Then

\vs{5pt}

\eqqed{
\phi(v)=0^{m_k-1} 1 0^{m_{k-1}-1} 1 \dots 0^{m_1-1} 1 0^{m_0} 1^{n_0-1} 0 1^{n_1-1} 0 \dots 1^{n_{k-1}-1} 0 1^{n_k}.
}
\ele

\ble[\cite{ss:mp}]
\label{cor2.4:ss:mp}
If $v\in\{0,1\}^*$ and $w=\phi(v)$, then

\vs{5pt}

\eqqed{\des v= d(w).}
\ele

Our first generating function concerns words $v$ that start with a $1$, if nonempty. We shall let $$
L=\{\epsilon\}\cup\{1u:u\in\{0,1\}^*\}
$$ 
be the set of such words, where $\epsilon$ is the empty word. The generating function we consider is
$$\sum_{v\in L} q^{\maj v}t^{\des v}x^{|v|}.$$
Applying Foata's bijection, the above lemmas, and equation~\ree{invtosumparts} shows that this sum is equal to
$$
\sum_{w\in\phi(L)} q^{|\lambda(w)|}t^{d(w)}x^{|w|}.$$

We claim that  a word $w$ is in $\phi(L)$ if and only if $\be(w)$ has all parts of size less than $k=d(w)$.  Keeping the notation of Lemma~\ref{prop2.2:ss:mp}, we see that $v\in L$ is equivalent to $m_0=0$.  So, by the same result, this will happen precisely when the subword of $w=\phi(v)$ following the $k$th one contains exactly $k$ zeros.  This means that the $k$th east step of $w$ will be the step which ends at the southeast corner of $D(w)$.  And this is equivalent to having part sizes in $\be(w)$ smaller than $d(w)$.
Since zero is allowed as a part in both $\be(w)$ and $\rho(w)$, standard generating function manipulations show that our desired sum is
$$
\sum_{k\geq 0} \frac{q^{k^2}t^k x^{2k}}{(x)_k (x)_{k+1}},$$
where
$$(x)_k = (1-x)(1-qx)\cdots(1-q^{k-1}x).$$
\bth
\label{maj:231,321}
We have 
$$M(231,321;x)=\sum_{k\geq 0} \frac{q^{k^2}t^k x^{2k}}{(x)_k (x)_{k+1}}.$$
\eth
\pf
From the preceding discussion, it suffices to show that there is a bijection $f$ from $\Av_n(231,321)$ to $L_n=\{v\in L:|v|=n\}$ which is Des preserving.
Recall from the proof of Proposition~\ref{inv2:form} that $\si\in\Av_n(231,321)$ if and only if $\sigma=\iota_j[\iotacheck_{n_1},\iotacheck_{n_2},\dots ,\iotacheck_{n_j}]$, where the first element of $\iotacheck_{n_i}$ is a left-right maximum $m_i$. Clearly, these left-right maxima are the only descents of $\sigma$, when not followed immediately by another left-right maximum.

Now we define $f(\sigma)=v=v_1v_2\cdots v_n$ where $v_i=1$ if $\sigma(i)$ is a left-right maximum, and $v_i=0$ if not. Since the first element of $\sigma$ is a left-right maximum,  $v_1=1$ and thus $f(\si)\in L_n$.  It also follows from the previous paragraph that $\Des\sigma=\Des f(\sigma)$.  This map is easily seen to be a bijection, and so we are done.\hqedm

The second generating function is a sum over the set of words which end in a $0$ if nonempty, denoted $R=\{u0:u\in\{0,1\}^*\}\cup\{\epsilon\}$:
$$\sum_{v\in R}q^{\maj v}t^{\des v}x^{|v|}.$$
Again, to find a formula for this sum we apply $\phi$.  Lemma~\ref{prop2.2:ss:mp} shows that in fact $\phi(R)=R$, so we obtain
$$\sum_{w\in R}q^{|\lambda(w)|}t^{d(w)}x^{|w|}.$$
Clearly having $w$ end in a zero is equivalent to $\rho(w)$ having all parts positive.
Therefore our desired sum is
$$
\sum_{k\geq 0} \frac{q^{k^2}t^kx^{2k}}{(x)_{k+1}(qx)_k}.$$

\bth
\label{maj:312,321}
We have
$$M(312,321;x) = \sum_{k\geq 0} \frac{q^{k^2}t^kx^{2k}}{(x)_{k+1}(qx)_k}.$$
\eth
\pf
As usual, it  suffices to show that we can find a Des-preserving bijection $f$ from $\Av_n(312,321)$ to $R_n$, the length $n$ words in $R$.
Notice that $R_{180}(\{231,321\})=\{312,321\}$.  Thus, from the proof of Theorem~\ref{maj:231,321}, it follows that $\sigma\in\Av_n(312,321)$  if and only if $\sigma=\iota_j[\iotahat_{n_1},\iotahat_{n_2},\dots ,\iotahat_{n_j}]$.  It follows that $\sigma$ is uniquely determined by its right-left minima, with the only descents of $\sigma$ being formed by these right-left minima, when not adjacent.

We create our word $f(\sigma)=v_1v_2\cdots v_n$ by letting $v_i=0$ if $\sigma(i)$ is a right-left minimum, and $v_i=1$ otherwise. Since $\sigma(n)$ is always a right-left minimum, we have $v_n=0$ and therefore $v\in R_n$. Also,  $\Des\sigma=\Des f(\si)$ by construction.  This process is clearly reversible, so $f$ is a bijection and we are done.\hqedm

\section{The inversion polynomial for triples}
\label{3I}

As before, we will ignore any $\Pi$ with $\{123,321\} \subset \Pi$. Since $\Pi\winv\Pi'$ forces $\inv\Pi=~\inv\Pi'$ for subsets of $\fS_3$, we compute $\inv\Pi$ for a representative of the equivalences given by equation~\ree{winv:eq} and, since they are all distinct, deduce that no others exist.  Thus we have the following result.

\begin{thm}
\label{inv:3-trivial}
We have

\bea
\left[123,132,231\right]_{\inv} &=& \big\{\{123,132,231\},\{123,132,312\},\{123,213,231\},\{123,213,312\}\big\},\\
\left[132,231,321\right]_{\inv} &=& \big\{\{132,231,321\},\{132,312,321\},\{213,231,321\},\{213,312,321\}\big\},\\
\left[132,213,231\right]_{\inv} &=& \big\{\{132,213,231\},\{132,213,312\}\big\},\\
\left[132,231,312\right]_{\inv} &=& \big\{\{132,231,312\},\{213,231,312\}\big\}.
\eea
All other  inv-Wilf equivalence classes for $\Pi\subset \fS_3$ with $\#\Pi=3$ and $\{123,321\}\not\sbs\Pi$ consist of a single triple.\hqedm
\end{thm}

As usual, in the above theorem, the results are listed so that each $[\Pi]_{\inv}$ is on a line adjacent to $[\Pi^{c}]_{\inv}$, and we only compute $I_n(\Pi;q)$ for one representative of each complementary pair.

\bpr
\label{inv:3-1}
For all $n\ge0$ we have 

\begin{align*}
&I_{n}(132,213,321;q)=\sum_{k=1}^{n}q^{k(n-k)},\\
&I_{n}(132,231,312;q)=\sum_{k=1}^{n}q^{k \choose 2},\\
&I_{n}(132,231,321;q)=[n]_{q},\\
&I_{n}(231,312,321;q)=\sum_{k=0}^{n}{n-k \choose k}q^{k}.
\end{align*}
\epr

\pf
The first three identities all follow from a straightforward count of inversions in the following characterizations of the permutations under consideration, where in each case $k$ varies over $[n]$: $\sigma\in\Av_{n}(132,213,321)$ if and only if $\sigma=21[\iota_{k},\iota_{n-k}]$; $\sigma\in\Av_{n}(132,231,312)$ if and only if $\sigma=12[\de_{k},\iota_{n-k}]$;  $\sigma\in\Av_{n}(132,231,321)$ if and only if $\sigma=213[1,\iota_{k-1},\iota_{n-k}]$.

For the fourth identity, we know from the proof of Proposition~\ref{inv2:form} that to avoid $231$ and $321$ we must have 
$\si=\iota_{j}[\iotacheck_{n_1},\iotacheck_{n_2},\dots,\iotacheck_{n_j}]$.  And to avoid $312$ we must have $n_i\le2$ for all $i$.  Thus $\sigma\in\Av_{n}(231,312,321)$ if and only if $\sigma=\iota_{j}[\de_{n_1},\de_{n_2},\ldots,\de_{n_j}]$ where $|\de_{n_i}|\le2$ for all $i$. With this characterization, the only inversions of $\sigma$ are descents within components of the form $\de_2$.  Now $\inv\sigma=k$ if and only if there are $k$ components of this form, and moreover $j=n-k$ since there are $k$ length two components and the rest are length one. 
So $\si$ is determined by choosing $k$ of the  $n-k$ components to be the ones of form $\de_2$.
The desired result follows.\hqedm
 
Consider the \emph{Fibonacci numbers} defined by $F_0=F_1=1$ and $F_n=F_{n-1}+F_{n-2}$ for $n\ge2$.  Simion and Schmidt~\cite{ss:rp} showed that $\#\Av_{n}(231,312,321)=F_n$.  So  the fourth equation in Proposition~\ref{inv:3-1} is a $q$-analogue of the well-known identity
$$
F_n=\sum_{k=0}^{n}{n-k \choose k}.
$$
A simple substitution of variables transforms $I_n(231,312,321)$ into the Fibonacci polynomials introduced by Morgan-Voyce~\cite{mor:lna}.  These have since been studied by many authors, with special attention given to their roots. 

\section{The major index polynomial for triples}
\label{3M}

Again, we will omit any $\Pi$ with $\{123,321\} \subset \Pi$. We have the following theorem for the remaining 16 triples and, as is our custom, the results are listed so that each $[\Pi]_{\maj}$ is adjacent to $[\Pi^{c}]_{\maj}$.
\begin{thm}
\label{maj:3-trivial}
We have

\bea
\left[123, 132, 312\right]_{\maj} &=& \big\{\{123, 132, 312\}\,\{123, 213, 231\},\{123, 231, 312\}\big\},\\
\left[132,213,321\right]_{\maj} &=& \big\{\{132,213,321\},\{132,312,321\},\{213,231,321\}\big\},\\
\left[132,213,231\right]_{\maj} &=& \big\{\{132,213,231\},\{132,231,312\}\big\},\\
\left[132,213,312\right]_{\maj} &=& \big\{\{132,213,312\},\{213,231,312\}\big\}.\\
\eea
All other  maj-Wilf equivalence classes for $\Pi\subset \fS_3$ with $\#\Pi=3$ and $\{123,321\}\not\sbs\Pi$ consist of a single triple.
\end{thm}
As before, it will suffice to prove the equivalences in the second and third equivalence classes above which we do  in the next proposition.

\bpr
\label{maj:triple}
For all $n\geq0$ and $\Pi_1= \{132,213,321\},$ $\{132,312,321\}$, or $\{213,231,321\}$ we have
$$
M_{n}(\Pi_1;q,t)= 1 + \sum_{k=1}^{n-1}q^{k}t = 1+qt[n-1]_q.
$$
For $\Pi_2= \{132,213,231\}$ or $\{132,231,312\}$ we have
$$
M_{n}(\Pi_2;q,t)= 1 + \sum_{k=2}^{n}{q^{k \choose 2}t^{k-1}}.
$$
We also have
\begin{align*}
&M_{n}(213,312,321;q,t)=1+(n-1)q^{n-1}t,\\
&M_{n}(132,231,321;q,t)=1+(n-1)qt.\\
\end{align*}
\epr
\pf 
The major index counts follow easily from the  characterizations of the  various avoidance classes given in the following table:  
$$
\barr{c|c}
\mbox{Class}		&\mbox{Elements}\\
\hline
\Av_{n}(132,213,321)	&21[\iota_{k},\iota_{n-k}]\\
\Av_{n}(132,312,321)	&213[\iota_{k},1,\iota_{n-k+1}]\\
\Av_{n}(213,231,321)	&132[\iota_{k-1},1,\iota_{n-k}]\\
\Av_{n}(132,213,231)	&21[\de_{k-1},\iota_{n-k+1}]\\
\Av_{n}(132,231,312)	&12[\de_{k},\iota_{n-k}]\\
\Av_{n}(213,312,321)	&132[\iota_{k-1},\iota_{n-k},1]\\
\Av_{n}(132,231,321)	&213[1,\iota_{k-1},\iota_{n-k}]
\earr
$$
where, as usual, $k$ ranges over $[n]$.
 \hqedm

For the last triple, there is no known closed form for its major index polynomial, so we will compute the generating function. Consider the set of words which are either the empty word or contain no consecutive ones and end in zero, denoted by
 $P$.
 The generating function we will investigate is
 $$\sum_{v\in P} q^{\maj v}t^{\des v}x^{|v|}.$$
Applying Foata's bijection, and using  Lemmas~\ref{thm4.3:foa:nin} and~\ref{prop2.2:ss:mp} as well as equation~\ree{invtosumparts}, we know the sum is equal to
$$
\sum_{w\in\phi(P)} q^{|\lambda(w)|}t^{d(w)}x^{|w|}.
$$
We claim that a word $w$ is in $\phi(P)$ if and only if $\rho(w)$, the portion to the right of  the Durfee square $D(w)$, is empty. This is equivalent to $w=\phi(v)$ ending in at least as many consecutive zeros as there  are total ones in $w$. Using the notation of Lemma~\ref{prop2.2:ss:mp}, we see that $v\in P$ if and only if $n_i = 1$ for $i<k$ and $n_k = 0$. The equation for $w=\phi(v)$ shows that there are at least $k$ zeros after the $k$th one since $1^{n_{i}-1}=1^0$. This gives us the desired property for $\phi(P)$. Also, zero is allowed as a part in $\be(w)$, the portion below $D(w)$, and thus our  sum becomes
$$
\sum_{k\ge0}\frac{q^{k^{2}}t^{k}x^{2k}}{(x)_{k+1}}.
$$
\bth
\label{maj:231,312,321}
We have 
$$
M(231,312,321;x)=\sum_{k\ge0}\frac{q^{k^{2}}t^{k}x^{2k}}{(x)_{k+1}}.
$$
\eth
\pf  From the above discussion, it suffices to show that there is a Des-preserving bijection $f$ from $\Av_{n}(231,312,321)$ to the set $P_n=\{v\in P:|v|=n\}$. 
Recall from the proof of Proposition~\ref{inv:3-1} that  $\sigma\in\Av_{n}(231,312,321)$ if and only if 
$\sigma=\iota_j[\de_{n_1},\de_{n_2},\dots ,\de_{n_j}],$ where $n_i \leq 2$ for all $i$.   Thus the right-left  minima of $\si$ are formed by the last elements of each $\de_{n_i}$.  And the only descents of $\sigma$ are formed by the first elements of the $\de_{n_i}$ where $n_i=2$ forcing no two to be adjacent.

We define $f(\sigma)=v=v_1v_2\cdots v_n$, where $v_i=0$ if $\sigma(i)$ is a right-left minimum and $v_i=1$ otherwise. Since the last element of $\sigma$ is always a right-left minimum, $v_n=0$. Also, there can be no consecutive ones in $v$ because of the description of the descents of $\si$ in the previous paragraph.  Hence, $f(\sigma)\in P_n$. By construction $\Des\sigma = \Des f(\sigma)$ and this process can be reversed to create $f^{-1}$.  This completes the proof.
\hqedm

It turns out that the polynomials $M_n(231,312,321;q)$ are precisely the $q$-Fibonacci numbers introduced by Carlitz~\cite{car:qfn} and subsequently studied by a number of authors.  In particular, they are closely related to pattern avoidance in set partitions as was first explored by Goyt and Sagan~\cite{gs:sps}

\section{Concluding remarks}
\label{cr}

\subsection{Other $\Pi\sbe\fS_3$}
\label{op}

We wish to complete the classification of the $\inv$- and $\maj$-Wilf equivalence classes for all $\Pi\sbe\fS_3$.  To do this, we must consider $\Pi$ with $\Pi\spe\{123,321\}$ or with $\#\Pi\ge4$.  

As already remarked, if  $\Pi\spe\{123,321\}$ then $\Av_n(\Pi)=\emp$ for $n\ge5$.   So these cases can be handled by mere computation.   It turns out that for certain $\Pi$ we also have $\Av_4(\Pi)=\emp$ and this leads to three sporadic $\maj$-Wilf equivalences.

If $\Pi\not\spe\{123,321\}$, then we use the following (slightly modified) result of Simion and Schmidt~\cite{ss:rp}.  In it, we let
$$
\tilde{\iota}_n=132[\iota_{n-2},1,1].
$$
\bpr[\cite{ss:rp}]
If $\{123,321\}\not\sbe \Pi\sbe\fS_3$ and  $\#\Pi=4$ then, for $n\ge2$, $\Av_n(\Pi)$ consists of exactly two elements taken from $\iota_n,\hat{\iota}_n,\tilde{\iota}_n$ and their rotations.\hqed
\epr
We note that $\maj R(\hat{\iota}_n)=\maj R(\tilde{\iota}_n)$ for any rotation $R\in D_4$.  This causes four more $\maj$-Wilf equivalences.

For $\#P=5,6$, at most $\iota_n$ or $\de_n$ survives in $\Av_n(\Pi)$, and so these sets are easy to classify.  Putting everything together, it is a simple (if tedious) matter to verify the following theorem.
\bth
Suppose $\Pi\sbe\fS_3$ and that either $\Pi\spe\{123,321\}$ or $\#\Pi\ge4$.
\ben
\item  All inv-Wilf equivalences among such $\Pi$ are trivial.
\item  We have the following maj-Wilf equivalence classes:
\bea
\left[123,132,213,231\right]_{\maj} &=& \big\{\{123,132,213,231\},\{123,132,231,312\}\big\},\\
\left[123,132,213,312\right]_{\maj} &=& \big\{\{123,132,213,312\},\{123,213,231,312\}\big\},\\
\left[123,132,312,321\right]_{\maj} &=& \big\{\{123,132,312,321\},\{123,213,231,321\}\big\},\\
\left[132,213,231,321\right]_{\maj} &=& \big\{\{132,213,231,321\},\{132,231,312,321\}\big\},\\
\left[132,213,312,321\right]_{\maj} &=& \big\{\{132,213,312,321\},\{213,231,312,321\}\big\},\\
\left[123,132,213,231,321\right]_{\maj} &=& \big\{\{123,132,213,231,321\},\{123,132,231,312,321\}\big\},\\
\left[123,132,213,312,321\right]_{\maj} &=& \big\{\{123,132,213,312,321\},\{123,213,231,312,321\}\big\},\\
\eea
and all other maj-Wilf equivalence classes consist of a single set.\hqedm
\een
\eth

\subsection{Mahonian Pairs}
\label{mp}

In~\cite{ss:mp}, Sagan and Savage introduced the notion of a Mahonian pair.  Let $S$ and $T$ be finite subsets of $\bbN^*$.  Then $(S,T)$ is a \emph{Mahonian pair} if the distribution of maj over $S$ is the same as the distribution of inv over $T$, equivalently,
$$
\sum_{s\in S}q^{\inv s}=\sum_{t\in T} q^{\inv t}
$$
So, for example, a special case of a celebrated theorem of MacMahon~\cite[pp.\ 508--549 and pp. 556--563]{mac:cp1} can be restated as saying that $(\fS_n,\fS_n)$ is a Mahonian pair and such pairs where maj and inv are equidistributed over $S=T$ have been the object of intense study.  But~\cite{ss:mp} is the first investigation of the case $S\neq T$.  However, all the pairs found in that paper are using words over $\{0,1\}^*$ and the authors asked if there were interesting pairs for larger alphabets.  We can now answer that question in the affirmative.  The following result follows immediately by comparing  the first two theorems in Section~\ref{2I} with the first two in Section~\ref{2M}, as well as the first two theorems in Section~\ref{3I} with the first three in Section~\ref{3M}.

\bth
\label{mp:thm}
Consider any $\Pi,\Pi'$ picked from corresponding rows in the same box of Table~\ref{mp:tab}.  For all $n\ge0$, if $S=\Av_n(\Pi)$ and $T=\Av_n(\Pi')$ then  $(S,T)$ is a Mahonian pair.\hqed
\eth

\begin{table}
$$
\barr{c|l}
\Pi & \{132,213\}, \{132,312\}, \{213,231\}, \{231,312\}\rule{0pt}{15pt}
\\ 
\Pi' &  \{132,231\}, \{132,312\}, \{213,231\}, \{213,312\}\rule{0pt}{15pt}
\\ \hline
\Pi & \{132,213,231\}, \{132,231,312\}\rule{0pt}{15pt}
\\ 
\Pi' & \{132,231,312\},\{213,231,312\}\rule{0pt}{15pt}
\\ \hline
\Pi & \{132,213,312\},\{213,231,312\}\rule{0pt}{15pt}
\\ 
\Pi' & \{132,213,231\},\{132,213,312\}\rule{0pt}{15pt}
\\ \hline
\Pi & \{123, 132, 312\},\{123, 213, 231\},\{123, 231, 312\}\rule{0pt}{15pt}
\\
\Pi' & \{123,132,231\},\{123,132,312\},\{123,213,231\},\{123,213,312\}\rule{0pt}{15pt}
\\ \hline
\Pi & \{132,213,321\},\{132,312,321\},\{213,231,321\}\rule{0pt}{15pt}
\\
\Pi' & \{132,231,321\},\{132,312,321\},\{213,231,321\},\{213,312,321\}\rule{0pt}{15pt}
\earr
$$
\caption{Mahonian pairs}
\label{mp:tab}
\end{table}

\subsection{Future work}

There are many directions in which this research could be taken.  One could  look at $\Pi$ containing patterns of length 4 or more.  These would give $q$-analogues of other famous combinatorial sequences which appear in pattern avoidance.  For example, it was proved by West~\cite{wes:gtc} that $\#\Av_n(3142,2413)=S_{n-1}$, a large Schr\"oder number.  These numbers are closely connected to the Catalan numbers and have many interesting properties.

One could investigate other statitics.  Just as inv and maj are equidistributed over $\fS_n$, so are another pair of well-known statistics.  One of them, des, we have already met.  The other is the number of excedences
$$
\exc\si=\#\{i\ |\ \si(i)>i\}.
$$
Statistics with the same distribution as des or exc are called \emph{Eulerian}.  Interesting resuts about the distribution of exc have been found by Elizalde~\cite{eli:fpe, eli:mpa}, Elizalde and Deutsch~\cite{ed:sub}, and Elizalde and Pak~\cite{ep:brr}.  And of course we already have results for des in the polynomials $M_n(\Pi;q,t)$.  But there are a host of other statistics which might give interesting results.  See the paper of Babson and Steingr{\'{\i}}msson~\cite{bs:gpp} for a partial list.

Rather than classical avoidance, one could look at generalized pattern avoidance.  A generalized pattern $\pi$ is one where, in order to have a copy of $\pi$ in $\si$ certain elements of the diagram of the copy must be adjacent either horizontally or vertically.  This concept was introduced in the paper of Babson and Steingr{\'{\i}}msson just cited.  They are of interest, among other reasons, because many statistics on permutations can be expressed as linear combinations of the avoidance statistics for these generalized patterns.

Finally, one could look at pattern avoidance in other combinatorial structures.  As previously mentioned, the $q$-Fibonacci numbers of Carlitz can also be obtained by avoiding patterns in set partitions.  Dokos, Dwyer, and Sagan are currently pursuing this line of research.

\medskip

\emph{Note added in proof.}  Recently Cheng, Elizalde, Kasraoui, and Sagan~\cite{ceks:imi} proved Conjecture~\ref{I321:con} and also answered Question~\ref{M321:ques} in the affirmative.  Also, Killpatrick~\cite{kil:wec} has proved Conjecture~\ref{M321parity}.

\medskip

\emph{Acknowledgement.}  We would like to thank two anonymous referees for carefully reading this article and suggesting useful changes.

\end{document}